\documentstyle[11pt]{article}
\begin{document}
\begin{center}
{\LARGE{\bf  Riemann Extension of the Anti- Mach  Space time  }}\\[1em]
\large{\bf{M. Abdel-Megied$^{*}$ \footnote{E-mail : amegied@frcu.eun.eg}, Nassar H. Abdel-All$^{**}$ \footnote{E-mail :nhabdeal2002@yahoo.com}\\ and E. A. Hegazy$^{*}$ \footnote{E-mail : elsayed.mahmoud@mu.edu.eg}}}\\
\normalsize {$^{*}$Mathematics Department, Faculty of Science, Minia University}\\
\normalsize {$^{**}$ Mathematics Department, Faculty of Science,
Assiut University}

\end{center}

\begin{abstract} Riemann extension for the anti Mach metric is
derived, the solution of geodesic equations for the extended space
are given, some properties for the extended space was studied and
compared with the basic space and the constructions of a translation
surface for the anti Mach metric in four dimension is established.

\end{abstract}

\setcounter{equation}{0}
\section{Introduction}
 Riemann extension is defined as a product of a non- Riemannian
$n$- space and a vector $n$-space $A^n$ (also called Riemannian $2n$
space) \cite{1}. It is possible by means of such extension, to
relate the properties of a non - Riemannian $n$ space with those of
certain Riemannian $2n$ spaces \cite{2}, \cite{z}.

For $n$ dimensional Riemannian space of the coordinates $x^i, \,\,
i=1:n$ the metric is given by:
\begin{equation}\label{d1}
    ^{n} ds^2= g_{ij} dx^i dx^j.
\end{equation}

 The Riemann extension of an affine connected space $A^n$ is given by the canonical form \cite{z}
\begin{equation}\label{r33}
   ^{2n} ds^2=  -2\Gamma_{ij}^{k} \Psi_{k} dx ^i dx ^ j + 2 d\Psi_{k}
    dx^{k},
\end{equation}
where $\Psi_k,\,\, k=n+1,...,2n$ are the coordinates of additional
space and $\Gamma_{ij}^{k}=\Gamma_{ji}^{k}$ are the connection
coefficients of the
space $A^n$.\\
 The geodesic equations of the metric (\ref{r33}) are given by:
\begin{equation}\label{r12}
    \ddot{y}^k+\hat{\Gamma}_{ij}^k \dot{y}^{i} \dot{y}^{j}=0,
\end{equation}
where $y^k$ are the coordinates of the extended space,
$\hat{\Gamma}_{ij}^k$ are the christoffel symbols of the second
kind for the metric (\ref{r33}) and dot means differentiation with respect to the affine parameter $s$.\\
 Dryuma \cite {5}-\cite{3} studied solutions
of the geodesic equations in the extended Riemanniann space and
gives some solutions of the Einstein equations for the eight
dimensional Riemann extension of the classical Schwarchild  four
dimensional metrics , Minkowsky space time metric in rotating
coordinate system, and some properties of the G\"{o}del space metric
and its Riemann extension .  We study the Riemann extension of the
Anti Mach metric space time and its geodesic equations in eight
dimension space with its complete solutions.  In \S{2} we introduce
the Anti Mach space time and its Riemann extension. The complete
solution of the geodesic equations in extended space time is given
in \S{3}. The constructions of translation surface for the Anti Mach
metric in four dimension is established in \S 4.

\setcounter{equation}{0}
\section{\bf{The Riemann Extension of the Anti-Mach Space time}}
The line element for the anti- Mach space time is given by
\cite{sch}:
\begin{equation}\label{r1}
   ^4 ds^2= dx^2- 4 t dx dz + 2 dy dz + 2 t^2 dz^2+ dt^2.
\end{equation}
 The non vanishing christoffel  symbols of second kind for the line
element (\ref{r1}) are given by:
\begin{equation}\label{r2}
    \Gamma_{13}^{4} = 1, \,\,\,\, \Gamma_{14}^{2} = -1,\,\,\,\, \Gamma_{33}^{4} = - 2
    t,\,\,\,\, \Gamma_{34}^{1} = -1.
\end{equation}
From equations (\ref{r1}), (\ref{r2}) and (\ref{r33}) the extended
metric takes the form:

\begin{equation}\label{r4}
\begin{array}{ccc}
    ^8 ds^2= 4P dzdt+ 4 Q dxdt - 4 Vdxdz + 4 tV dz^2 \\+ 2 dxdP + 2
    dydQ+ 2 dzdU+ 2 dt dV\\
    =\hat{g}_{ij} dy^i dy^j
    \end{array}
,\end{equation}
where $\Psi_k = (P,Q,U,V) $ are the additional coordinates.\\
  In the local coordinates $y^i=(x,y,z,t,P,Q,U,V)$ of  this type
of metrics are metrics with vanishing curvature invariant
\cite{sch}, \cite{Ram} . They play an important role in general
theory of relativity, in particular for the pp-waves \cite{JO}. The
eight dimensional space (\ref{r4}) is also Einstein space with the
condition that the Ricci tensor is vanishing . i.e
\begin{equation}\label{xc22}
 ^8R_{ik}=0
\end{equation}

From equation (\ref{r12}) the complete system of geodesic equations
for the metric (\ref{r4})  read as:
\begin{equation}\label{rr4}
   \ddot{x}= 2 \dot{z}\dot{t},
\end{equation}
\begin{equation}\label{r5}
 \ddot{y}= 2 \dot{x}\dot{t},
\end{equation}
\begin{equation}\label{r6}
     \ddot{z}=0
\end{equation}
\begin{equation}\label{r7}
    \ddot{t}= 2t \dot{z}^2- 2 \dot{x} \dot{z},
\end{equation}
\begin{equation}\label{r8}
    \ddot{P}= 4 Q \dot{z} (\dot{x}- t \dot{z})- 2 \dot{t} \dot{Q}+ 2
    \dot{z} \dot{V}
\end{equation}
\begin{equation}\label{r9}
    \ddot{Q}=0
\end{equation}
\begin{equation}\label{r10}
    \ddot{U} = 4 P \dot{z} (\dot{x}- t \dot{z}) - 2 \dot{t} \dot{P}+
    2(\dot{x}- 2t \dot{z}) \dot{V}
\end{equation}
\begin{equation}\label{r11}
    \ddot{V}= 2V \dot{z}^2 - 4 Q \dot{z} \dot{t} -2\dot{z}\dot{P}- 2
    \dot{x}\dot{Q},
\end{equation}
where dot denotes differentiation with respect to the parameter
$s$.\\
 The four equations (\ref{rr4}) - (\ref{r7})
constitute the system of geodesics equations of the basic space and
not contain the coordinates $\Psi_k$ while the other four equations
(\ref{r8})- (\ref{r11})  represent a system of second order ordinary
differential equations for the coordinates $\Psi_{k}$ and the first
derivative of the coordinates $(x,y,z,t)$.
 Solutions of the equations (\ref{rr4})- (\ref{r11}) will be
given in the next section.
 \setcounter{equation}{0}
\section {\bf{Solutions of the Geodesics Equations}}
To solve the geodesic equations (\ref{rr4}) - (\ref{r11}) we take an
arbitrary initial point $x^{i}(s)|_{s=0}=x^{i}_{0}$ and
$\Psi_{k}(s)|_{s=0}=\Psi_{0}=(P_0,Q_0,U_0,V_0)$ with the initial
directions $\frac{dx^{i}}{ds}|_{s=0} = \xi^{i}$ and
$\frac{d\Psi_{k}}{ds}|_{s=0} = h_{k}$, with the supplementary
condition for the  geodesics
\begin{equation}\label{rr42}
  \hat{g}_{ij}\frac{dy^i}{ds}\frac{dy^j}{ds}=\epsilon,
\end{equation}
where $\epsilon= 0 , 1, -1$ according to null, timelike and
spacelike geodesics respectively, that is
\begin{equation}\label{d1}
    4P \dot{z}\dot{t}+ 4 Q \dot{x} \dot{t} - 4 V \dot{x} \dot{z} + 4 tV \dot{z}^2 + 2 \dot{x} \dot{P} + 2
    \dot{y} \dot{Q}+ 2 \dot{z} \dot{U}+ 2 \dot{t} \dot{V}=\epsilon
.\end{equation} Solutions of the first four equations (\ref{rr4})-
(\ref{r7}) are given by Ozsv\'{a}th and Sch\"{u}cking  \cite{sch} and read as:\\
For $\xi^3 \neq 0$
\[ x=x_0+ (2 t_0 \xi^{3}- \xi^{1}) s + \frac{\xi^4}{\xi^3} (1- \cos(\sqrt{2}
    \xi^{3})s)\]\begin{equation}\label{r20}+ \frac{\sqrt{2}(\xi^{1}- t_0 \xi^3)}{\xi^3} \sin(\sqrt{2}
    \xi^{3})s
\end{equation}
\[y= y_0 + [\xi^2 + \frac{(\xi^4)^2+ 2 (\xi^1)^2}{2 \xi^3}+ t_0
(3t_0\xi^{3}-4\xi^1)]s +\]\[\,\,\,\,\,\,\, \frac{\xi^4}{2
    (\xi^3)^2}[2(\xi^1-2 t_0 \xi^3)\cos\sqrt{2}\xi^3 s- (\xi^1- t_0 \xi^3) \cos 2\sqrt{2} \xi^3 s\]\[-(\xi^1-3 t_0
    \xi^3)]-\frac{\sqrt{2}}{(\xi^3)^2}[(\xi^1)^2-t_0 \xi^3 (3\xi^1-2t_0
    \xi^3)] \sin\sqrt{2}\xi^3 s\]\begin{equation}\label{r21}+\frac{1}{2\sqrt{2}(\xi^3)^2}[(\xi^1)^2- \frac{(\xi^4)^2}{2}- t_0 \xi^3(2\xi^1-
    t_0\xi^3)]\sin2\sqrt{2}\xi^3 s
\end{equation}
\begin{equation}\label{r22}
    z= z_0 + \xi^{3} s
\end{equation}
\begin{equation}\label{r23}
    t= t_0 (2- \cos \sqrt{2}\xi^3 s)- \frac{\xi^1}{\xi^3}(1- \cos\sqrt{2}\xi^3
    s)+ \frac{\xi^4}{\sqrt{2}\xi^3} \sin\sqrt{2}\xi^3 s,
\end{equation}
and for $\xi^3 = 0$
\begin{equation}\label{v1}
\begin{array}{cc}
x=x_0+\xi^1 s\\
y=y_0 +\xi^2 s +\xi^1 \xi^4 s^2\\
z=z_0\\
t=t_0+\xi^4 s
\end{array}
\end{equation}
    To solve equations (\ref{r8})-
(\ref{r11}) for the case ($\xi^3\neq0$) we have from equation
(\ref{rr4})
\begin{equation}\label{r24}
\dot{x}(s)= 2 \xi^3 t(s) + \xi^1- 2\xi^3t_0.
\end{equation}
Solution of the equation (\ref{r9}) is:
\begin{equation}\label{r25}
    Q= h_2 s+ Q_0.
\end{equation}
Using equations (\ref{r22}), (\ref{r24}) and (\ref{r25}) in equation
(\ref{r8}) we have:
\[
    \ddot{P}= 4Q\xi^{3}[(\xi^1-\xi^3 t_0) \cos (\sqrt{2} \xi^3) s +\]\[ \frac{\xi^4}{\sqrt{2}} \sin(\sqrt{2} \xi^{3}) s]
    + 4\xi^3 h_2 (\xi^1-\xi^3 t_0) s \cos(\sqrt{2} \xi^3) s\]\begin{equation}\label{r26}+
    2\sqrt{2} \xi^3 \xi^4 h_2 s \sin(\sqrt{2} \xi^3)s-2h_2 \dot{t}+
    2\xi^3 \dot{V}
,\end{equation}  Integrating with respect to $s$ we get:

\[
    \dot{P}= 2\sqrt{2} Q_0 (\xi^1-\xi^3t_0) \sin(\sqrt{2}\xi^3)s -
    2Q_0\xi^4\cos(\sqrt{2}\xi^3)s+\]\[4\xi^3 h_2(\xi^1-\xi^3 t_0)[\frac{s}{\sqrt{2}\xi^3} \sin(\sqrt{2}\xi^3)s+ \frac{1}{2(\xi^3)^2}\cos(\sqrt{2} \xi^3)s]
+ \]\[2\sqrt{2} \xi^3 \xi^4 h_2[\frac{-s}{\sqrt{2} \xi^3} \cos
(\sqrt{2}\xi^3)s +\frac{1}{2 (\xi^3)^2}
\sin(\sqrt{2}\xi^3)s]\]\begin{equation}\label{r27}- 2 h_2 t(s)+2
\xi^3 V+ C_1 ,\end{equation} where $C_1$ is a constant of
integration. Using equation (\ref{r23}) in equation (\ref{r27}) we
get:

\[
    \dot{P}= L_1 \sin(\sqrt{2}\xi^3)s+ L_2 \cos(\sqrt{2}\xi^3)s +
    L_3 s \sin(\sqrt{2}\xi^3)s+\]\begin{equation}\label{r28}L_4 s \cos(\sqrt{2}\xi^3)s+ L_5 +
    2\xi^3 V+ C_1,
\end{equation}
where
$$L_1 = 2\sqrt{2} Q_0 (\xi^1-\xi^3 t_0),\quad L_2= -2 Q_0 \xi^4,\quad L_3 =2\sqrt{2} h_2(\xi^1-\xi^3 t_0),$$ $$ L_4= -2\xi^4
h_2,\quad L_5 = \frac{2h_2\xi^1}{\xi^3}-4h_2t_0.$$ From the initial
conditions (at $s\rightarrow0 $ we have $V=V_0$ and $\dot{P}= h_1$)
in equation (\ref{r28}) we get:
\[ L_5+ C_1 = h_1-L_2-2\xi^3 V_0,\]
then we have:
\[
    \dot{P}= L_1 \sin(\sqrt{2}\xi^3)s+ L_2 (\cos(\sqrt{2}\xi^3)s-1) +
    L_3 s \sin(\sqrt{2}\xi^3)s+\]\begin{equation}\label{r29}L_4 s \cos(\sqrt{2}\xi^3)s+
    2\xi^3 (V-V_0)+ h_1,\end{equation}

 using equations (\ref{r22}), (\ref{r23}), (\ref{r25}) and (\ref{r29}) in equation (\ref{r11})
we get:
\begin{equation}\label{r30}
    \ddot{V}+2 (\xi^3)^2 V= M+M_1 \sin(\sqrt{2}\xi^3)s +M_2
    \cos(\sqrt{2}\xi^3)s ,
\end{equation}
where $$M=2 \xi^1 h_2 - 2 \xi^3 (h_1 + 2 \xi^4 Q_0 + 2 h_2 t_0) + 4
(\xi^3)^2 V_0,$$ $$M_1= -2\sqrt{2}\xi^4h_2, \quad M_2 = -4h_2
(\xi^1-\xi^3 t_0).$$
 The solution of equation (\ref{r30}) is given by:

   \[ V(s) = A_1\sin(\sqrt{2} \xi^3)s+A_2\cos(\sqrt{2} \xi^3)s+\frac{M}{2
   (\xi^3)^2}\]
\begin{equation}\label{r31}
    -\frac{M_1}{2\sqrt{2} \xi^3} s \cos(\sqrt{2}\xi^3)s+\frac{M_2}{2\sqrt{2} \xi^3}  s \sin(\sqrt{2}\xi^3)s
,\end{equation} where $ A_1$ and $A_2$ are constants.\\
From the initial condition  $s=0$ we have:
$$V=V_0, \quad \dot{V}= h_4 ,$$
and from the  equation (\ref{r31}) we get:
\[A_1= \frac{h_4}{\sqrt{2}\xi^3}+\frac{M_1}{4 (\xi^3)^2},\quad A_2= V_0 -\frac{M}{2(\xi^3)^2} \]
Hence equation (\ref{r31}) becomes :
\[
    V(s) = K_1 \sin(\sqrt{2}\xi^3)s
    +K_2 \cos(\sqrt{2}\xi^3)s +K_3 s \sin(\sqrt{2}\xi^3)s
    \]\begin{equation}\label{r32}+K_4 s \cos(\sqrt{2}\xi^3)s ,
\end{equation}where
$$K_1= \frac{h_4}{\sqrt{2}\xi^3}+\frac{M_1}{4 (\xi^3)^2},\quad K_2=V_0 -\frac{M}{2(\xi^3)^2}, \quad K_3= \frac{M_2}{2\sqrt{2} \xi^3},\quad K_4=- \frac{M_1}{2\sqrt{2} \xi^3}.$$
Using the equation (\ref{r32}) in the equation (\ref{r29}) we get:

\begin{equation}
    \dot{P}= H_1 \sin(\sqrt{2} \xi^3)s+H_2 \cos(\sqrt{2}
    \xi^3)s+H_3
,\end{equation} where
\[ H_1= \frac{1}{2 \sqrt{2} \xi^3}[-\xi^4 h_2 +4 \xi^3 (h_4+2 Q_0 (\xi^1-\xi^3 t_0))], \quad H_2= 2(h_1 - \frac{\xi^1 h_2}{\xi^3}+\xi^4 Q_0 + 2 h_2 t_0-\xi^3V_0)\]\[ H_3= -h_1-2\xi^4 Q_0+  \frac{2 \xi^1 h_2}{\xi^3}-4 h_2 t_0+2 \xi^3 V_0.\]
 By integration with respect to $s$ we obtain:

     \begin{equation}\label{r344}P= \frac{- H_1}{\sqrt{2} h_3} \cos(\sqrt{2} h_3)s+\frac{H_2}{\sqrt{2} h_3} \sin(\sqrt{2}
    h_3)s+ H_3 s+ C_2
,\end{equation} where $ C_2$ is a constant.
 From the initial condition  at $s=0$ we
have $P=P_0$ then:
\[C_2 = P_0+ \frac{H_1}{\sqrt{2}\xi^2}.\]
Equation (\ref{r344}) read as:

     \begin{equation}\label{r34}P= P_0+\frac{H_1}{\sqrt{2} \xi^3}(1- \cos(\sqrt{2} \xi^3)s)+\frac{H_2}{\sqrt{2} \xi^3} \sin(\sqrt{2}
    \xi^3)s+ H_3 s.
\end{equation}
In the following we integrate the equation for the coordinate $U$,
  using the equations
(\ref{r20}), (\ref{r22}), (\ref{r23})
 , (\ref{r24}), (\ref{r32}) and (\ref{r34}) in the equation (\ref{r10}) we get:

\[  \ddot{U}= R_1 \sin(\sqrt{2}\xi^3)s +R_2 \cos(\sqrt{2}\xi^3)s
    +R_3 s \sin (\sqrt{2}\xi^3)s\]\[ + R_4 s \cos(\sqrt{2}\xi^3)s
    +R_5 \sin^2(\sqrt{2} \xi^3)s\]\begin{equation}\label{r36}+R_6 \cos^2(\sqrt{2} \xi^3)s+R_7 \sin(\sqrt{2} \xi^3)s \cos(\sqrt{2}
    \xi^3)s,\end{equation}
where $R_1$, $R_2$, $R_3$, $R_4$, $R_5$, $R_6$,
 and $R_{7}$ are constants read as:
 \[R_1=\frac{1}{\sqrt{2} \xi^3}[8 (\xi^1)^2 h_2-(\xi^4)^2 h_2+4 \xi^3 \xi^4 (h_4+\xi^3 (P_0+4 Q_0 t_0))+4 (\xi^3)^2 t_0 (3 h_1+6
h_2 t_0\]\[-4 \xi^3 V_0)+4 \xi^1 \xi^3 (-2 h_1-2 \xi^4 Q_0-7 h_2
t_0+3 \xi^3 V_0)] ,\]
\[R_2=\frac{1}{2 \xi^3}[-7 \xi^1 \xi^4 h_2-8 (\xi^3)^3 t_0 (P_0-2 Q_0 t_0)+4 \xi^3 (2 (\xi^4)^2 Q_0+\xi^1 (3 h_4+4 \xi^1 Q_0)+\xi^4
(h_1+3 h_2 t_0))\]\[+8 (\xi^3)^2 (-2 h_4 t_0+\xi^1 (P_0-4 Q_0
t_0)-\xi^4 V_0)]
 ,\]
\[R_3= -2 \sqrt{2} \xi^4 (-\xi^1 h_2+\xi^3 (h_1+2 \xi^4 Q_0+2 h_2 t_0)-2 (\xi^3)^2 V_0)
,
\]
\[R_4=4 (\xi^1-\xi^3 t_0) (\xi^1 h_2-\xi^3 (h_1+2 \xi^4 Q_0+2 h_2 t_0)+2 (\xi^3)^2 V_0)
\]
\[R_5= \frac{1}{\xi^3}[-5 \xi^1 \xi^4 h_2+8 (\xi^3)^3 Q_0 t_0^2+\xi^3 (4 (\xi^4)^2 Q_0+4 \xi^1 (h_4+2 \xi^1 Q_0)+\xi^4 (4 h_1+9 h_2 t_0))\]\[-4
(\xi^3)^2 (h_4 t_0+4 \xi^1 Q_0 t_0+\xi^4 V_0)]
 ,\]
\[R_6=\frac{1}{\xi^3}[5 \xi^1 \xi^4 h_2-8 (\xi^3)^3 Q_0 t_0^2-\xi^3 (4 (\xi^4)^2 Q_0+4 \xi^1 (h_4+2 \xi^1 Q_0)\]\[+\xi^4 (4 h_1+9 h_2 t_0))+4
(\xi^3)^2 (h_4 t_0+4 \xi^1 Q_0 t_0+\xi^4 V_0)], \]
\[R_7=\frac{\sqrt{2}}{\xi^3}[ (-8 \xi^1)^2 h_2+(\xi^4)^2 h_2-4 \xi^3 \xi^4 h_4+8 \xi^1 \xi^3 (h_1+3 h_2 t_0-\xi^3 V_0)+8 (\xi^3)^2 t_0 (-h_1-2
h_2 t_0+\xi^3 V_0)] .\]

  Integrating twice and using the initial conditions we get:
 \[U(s)= U_0+ N_1 \sin(\sqrt{2}\xi^3)s++N_2 (\cos(\sqrt{2}\xi^3)s-1)
    \]\[\quad\quad\quad\quad\quad\quad+ N_3 s \sin(\sqrt{2} \xi^3)s+N_4 s \cos(\sqrt{2} \xi^3)s+ N_5 \sin(2\sqrt{2}
    \xi^3)s\]\begin{equation}\label{r37}+N_6 ( \cos(2\sqrt{2} \xi^3)s-1)
    + N_7 s,
\end{equation}
where $N_1$, $N_2$, $N_3$, $N_4$, $N_5$, $N_6$ and $N_7$ are
constants and are given by:
\[N_1=\frac{1}{2 \sqrt{2} (\xi^3)^3}[(\xi^4)^2 h_2-4 \xi^3 \xi^4 (h_4+\xi^3 P_0+2 \xi^1 Q_0)+4 \xi^3 (\xi^1 h_2 t_0+\xi^3 (-h_1 t_0-2
h_2 t_0^2+\xi^1 V_0))]
 ,\]
\[N_2=\frac{1}{4 (\xi^3)^3}[-\xi^1 \xi^4 h_2+8 (\xi^3)^3 t_0 (P_0-2 Q_0 t_0)+4 \xi^3 (2 (\xi^4)^2 Q_0-\xi^1 (3 h_4+4 \xi^1 Q_0)\]\[+\xi^4
(h_1+h_2 t_0))-8 (\xi^3)^2 (-2 h_4 t_0+\xi^1 (P_0-4 Q_0 t_0)+\xi^4
V_0)]
 ,\]
\[N_3=\frac{\sqrt{2} \xi^4 (-\xi^1 h_2+\xi^3 (h_1+2 \xi^4 Q_0+2 h_2 t_0)-2 (\xi^3)^2 V_0)}{(\xi^3)^2}
,\]

\[N_4=-\frac{2 (\xi^1-\xi^3 t_0) (\xi^1 h_2-\xi^3 (h_1+2 \xi^4 Q_0+2 h_2
t_0)+2 (\xi^3)^2 V_0)}{(\xi^3)^2},\]
\[N_5= \frac{1}{8 \sqrt{2} (\xi^3)^3}[8 (\xi^1)^2 h_2-(\xi^4)^2 h_2+4 \xi^3 \xi^4 h_4+8 (\xi^3)^2 t_0 (h_1+2 h_2 t_0-\xi^3 V_0)\]\[+8 \xi^1 \xi^3 (-h_1-3
h_2 t_0+\xi^3 V_0)] ,\]
\[N_6=\frac{1}{8 (\xi^3)^3}[-5 \xi^1 \xi^4 h_2+8 (\xi^3)^3 Q_0 t_0^2+\xi^3 (4 (\xi^4)^2 Q_0+4 \xi^1 (h_4+2 \xi^1 Q_0)+\xi^4 (4 h_1+9 h_2 t_0))\]\[-4
(\xi^3)^2 (h_4 t_0+4 \xi^1 Q_0 t_0+\xi^4 V_0)] ,\] and
\[N_7=-\frac{(\xi^4)^2 h_2}{4 (\xi^3)^2}+h_3+2 \xi^4 P_0+2 h_1 t_0+4 \xi^4 Q_0 t_0+4 h_2 t_0^2+\frac{\xi^4 h_4-2 \xi^1 h_2 t_0}{\xi^3}-2 \xi^3 t_0 V_0
,\]
 For the second case ($\xi^3=0$) it is easy to find:
\begin{equation}\label{v5}
\begin{array}{cc}
    P=P_0+h_1 s- \xi^4 h_2 s^2\\
 Q=Q_0 + h_2 s\\
  U=U_0-(\xi^4 h_1-\xi^1h_4) s^2+ \frac{2}{3}h_2[(\xi^4)^2 -
    (\xi^1)^2] s^3 + h_3 s\\
    V=V_0+h_4 s-\xi^1 h_2 s^2.
    \end{array}
\end{equation}
 These solutions for $\xi^3 \neq 0$ and $\xi^3 = 0$ satisfy also  the
 condition given by equation
 (\ref{d1}).\\
 Because of the homogeneity of the basic manifold ( also extended
  manifold) without loss of generality, we consider the geodesics  in eight dimension space time
  at the vertex (i.e $x_{0}^{i}=0, \Psi_0=0 )$.

For ($\xi^3\neq0$) the  equations  (\ref{r25}), (\ref{r32}),
(\ref{r34}) and (\ref{r37}) are reduced to:
\[
P=\frac{1}{4 (\xi^3)^2}[-4 s (\xi^3)^2 h_1+(-4 \sqrt{2}\xi^1 \sin
\sqrt{2} s \xi^3 +\xi^4(-1+\cos \sqrt{2} s \xi^3))
h_2+\]\begin{equation}\label{r41}4 \xi^3 (\sqrt{2} h_1 \sin \sqrt{2}
s \xi^3 +2 s \xi^1 h_2+h_4-h_4 \cos \sqrt{2} s \xi^3)]
\end{equation}
\begin{equation}\label{r42}
    Q=h_2 s
\end{equation}
\[U=
\frac{1}{16 (\xi^3)^3}[\sqrt{2}  (-8 \xi^1 \xi^3 h_1+8 (\xi^1)^2
h_2-(\xi^4)^2 h_2+4 \xi^3 \xi^4 h_4) \sin 2 \sqrt{2} s \xi^3\]\[-4
\sqrt{2}  (-4 s (\xi^3)^2 \xi^4 h_1-(\xi^4)^2 h_2+4 \xi^3 \xi^4 (s
\xi^1 h_2+h_4))\sin \sqrt{2} s \xi^3+\]\[2 (7 \xi^1 \xi^4 h_2+8 s
(\xi^3)^3 h_3+8 s (\xi^3)^2 \xi^4 h_4-\xi^3 (\xi^4 (12 h_1+2 s \xi^4
h_2)-20 \xi^1 h_4))\]\[+2 (-5 \xi^1 \xi^4 h_2+\xi^3 (4 \xi^4 h_1+4
\xi^1 h_4)) \cos 2 \sqrt{2} s\xi^3\]\begin{equation}\label{r43}-4
 (-4 \xi^3 \xi^4 h_1+8 s (\xi^1)^2 \xi^3
h_2+\xi^1 (\xi^4 h_2+4 \xi^3 (-2 s \xi^3 h_1+3 h_4))) \cos \sqrt{2}
s \xi^3]
\end{equation}

\[
V(s)=\frac{1}{8 (\xi^3)^2}[-8 \xi^3 h_1+8 \xi^1 h_2+8  (-\xi^1
h_2+\xi^3 (h_1+s \xi^4 h_2))\cos \sqrt{2} s
\xi^3\]\begin{equation}\label{r44}+\sqrt{2} (-8 s \xi^1 \xi^3
h_2-\xi^4 h_2+4 \xi^3 h_4) \sin\sqrt{2} s \xi^3)].
\end{equation}
Equation (\ref{d1}) lead to:
\begin{equation}\label{v12}
   \hat{g}_{ij} \frac{dy^i}{ds}\frac{dy^j}{ds}|_{y^i=0}=2 \xi^1 h_1+2 \xi^2 h_2+\frac{3 (\xi^4)^2 h_2}{2 \xi^3}+2 \xi^3 h_3+2 \xi^4 h_4
=\epsilon.
\end{equation}

For $\xi^3=0$ the equation (\ref{d1}) take the form:
\begin{equation}\label{z1}
     \xi^1 h_1+\xi^2 h_2+ \xi^4 h_4=\epsilon,
\end{equation} also the equation (\ref{v5}) reduced to
\begin{equation}\label{vvv1}
\begin{array}{cc}
    P  =  h_1 s- \xi^4 h_2 s^2,\\
 Q =  h_2 s,\\
    U  =  -(\xi^4 h_1-\xi^1h_4) s^2+ \frac{2}{3}h_2[(\xi^4)^2 -
    (\xi^1)^2] s^3 + h_3 s,\\
    V  =  h_4 s-\xi^1 h_2 s^2.
 \end{array}
 \end{equation}
 The geometrical and physical properties of the geodesic
 equation, and the relation between the extended space and the basic space  will be
 done in the future work.
\setcounter{equation}{0}
\section{Translation Surfaces}
Translation surfaces in an arbitrary Riemannian space are defined by
the systems of equations \cite{Eis}, \cite{HU}
\begin{equation}\label{bb12}
    \frac{\partial^2 x^i(u,v)}{\partial u \partial v} + \Gamma^i _
   { jk} \frac{\partial  x^j(u,v)}{\partial u} \frac{\partial  x^k(u,v)}{\partial
    v}=0,
\end{equation}
where $\Gamma^i _{jk}$ is the christoffel symbols of the second kind.\\
    In the following we study the translation surfaces of the Anti-
    Mach space time metric given by the equation (\ref{r1}).
Using equation (\ref{r1}) in equation (\ref{bb12}) we get:
\begin{equation}\label{bb1}
    \frac{\partial^2 x}{\partial u \partial v}- \frac{\partial z}{\partial
    u} \frac{\partial t}{ \partial v}- \frac{\partial t}{\partial u}
    \frac{\partial z}{\partial v}=0
\end{equation}
\begin{equation}\label{bb2}
\frac{\partial^2 y}{\partial u \partial v}- \frac{\partial
x}{\partial
    u} \frac{\partial t}{ \partial v}- \frac{\partial t}{\partial u}
    \frac{\partial x}{\partial v}=0
\end{equation}
\begin{equation}\label{bb3}
    \frac{\partial^2 z}{\partial u \partial v}=0
\end{equation}
\begin{equation}\label{bb4}
    \frac{\partial^2 t}{\partial u \partial v}+ \frac{\partial x}{\partial
    u} \frac{\partial z}{ \partial v}+ \frac{\partial z}{\partial u}
    \frac{\partial x}{\partial v}- 2 t \frac{\partial z}{ \partial u}\frac{\partial z}{ \partial
    v}=0.
\end{equation}
The solution of the equation (\ref{bb3}) takes the form
\begin{equation}\label{bb5}
     z(u,v)= f(u)+ g(v).
\end{equation}
 Taking the function $t(u,v)$ as in the following
\begin{equation}\label{bb6}
    t(u,v)= f(u)- g(v).
\end{equation}
Using  equations (\ref{bb5}) and (\ref{bb6}) in equation (\ref{bb1})
we have the representation
\begin{equation}\label{bb7}
    x(u,v)= x_{1} (u) + x_{2} (v).
\end{equation}

From the equations (\ref{bb7}) and (\ref{bb1}) we get:
\begin{equation}\label{bb17}
    \frac{x_{1} '- 2 ff'}{f'}  = \frac{x_{2} '+ 2 gg'}{g'}=C_3,
\end{equation}
where $C_3$ is a constant, $f'= \frac{df(u)}{du}$ and $g'=
\frac{dg(v)}{dv}$.

By integration we get
\begin{equation}\label{bb9}
    x_1=  f^2+C_3 f+C_4, \quad\quad x_2 =  - g^2 - C_3 g+C_5,
\end{equation}
where $C_4$ and $C_5$ are constants of integration.

 From the equations (\ref{bb6}), (\ref{bb7}) and (\ref{bb9}) in the equation
(\ref{bb2}) we get:
\begin{equation}\label{bb10}
    \frac{\partial^2 y}{\partial u \partial v}= -2(C_3 f' g'+ ff'g'+
    gg'f').
\end{equation}
Solution of the equation (\ref{bb10}) is given by:
\begin{equation}\label{bb11}
    y(u,v)= -2 C_3 fg -f^2 g-g^2 f +G_1(u)+ G_2 (v)
\end{equation}
where $G_1(u)$ and $G_2 (v)$ are two arbitrary functions. It is
obvious from equation (\ref{bb11}) that $y$ is not a sum of two
functions, function depended on the coordinate $u$ and the other
depended on the coordinate $v$. Hence, there is no existence of a
translation surface in four dimension \cite{HU}, but its projection
in the 3- dimensional space $(x,z,t)$ is a translation surface and
we have the following theorem.

 \textbf{ Theorem:} \textit{For the anti-Mach metric: There
exists a translation surface in the orthogonal projection in the
3-dimensional space  $(x,z,t)$ and the existence comprises up to two
an arbitrary functions in one variable.}
\section{\bf{Conclusion}} We
have offered a Riemann extension for the anti-Mach space time. It is
found that the vacuum field equations in eight dimensional space
time are also satisfied . The full  geodesic equations in this
extended space time are obtained. We studding also the translation
surfaces for the basic space.

\end{document}